# OPTIMIZATION IN ENGINE DESIGN VIA FORMAL CONCEPT ANALYSIS USING NEGATIVE ATTRIBUTES

Rodríguez-Jiménez, J. M. [1], Cordero, P. [1], Enciso, M. [1] and Mora, A [1].

[1]Universidad de Málaga, Andalucía Tech, Boulevar Louis Pasteur SN, Málaga, Spain
jmrodriguez@ctima.uma.es

## ABSTRACT

*There is an exhaustive study around the area of engine design that covers different methods that try to reduce costs of production and to optimize the performance of these engines. Mathematical methods based in statistics, self-organized maps and neural networks reach the best results in these designs but there exists the problem that configuration of these methods is not an easy work due the high number of parameters that have to be measured.*

*In this work we extend an algorithm for computing implications between attributes with positive and negative values for obtaining the mixed concepts lattice and also we propose a theoretical method based in these results for engine simulators adjusting specific and different elements for obtaining optimal engine configurations.*

## 1. INTRODUCTION

An engine, or motor, is a machine designed to convert one form of energy into mechanical energy. There are several kinds of engines that are used in different situations. The most popular are heat engines, including internal combustion engines and external combustion engines, that burn fuel to create heat, which then creates a force. Some of these engines create electricity that is used in electric motors that convert electrical energy into mechanical motion. Engines are not only artificial machines because in nature we can find some examples in biological systems with the form of molecular motors, like myosin in muscles, that use chemical energy to create forces and eventually motion in animals and plants.

The performance of an engine measure different kind of properties like engine speed, torque, power, efficiency or sound levels. Focusing in have the best result in one of them could do that other properties have a poor result so engineers have to take into account all of these properties in general.

The cost of engine design is high so engineers have to use simulators to configure models. In the design of engines, simulators are an excellent tool to save cost of production because in them engineers can adjust different parameters and observe the global result. Complex engines could have thousands of parameters that have to be adjusted for obtaining optimal results in different ways. The problem is: what parameter they have to adjust? How these parameters have to be modified?

In general usage, design of experiments (DOE) or experimental design is the design of any information-gathering exercises where variation is present. Engineers could manage this process under full control, where statistics are usually used, or not, using artificial intelligence methods. Formal planned experimentation is often used in evaluating physical objects, structures, components and materials. Design of experiments is thus a discipline that has very broad application across all the natural and social sciences and, of course, engineering.

Datasets are necessary to collect all the information in experiments results. There are several kinds of datasets but in engine design all data is measured in numerical values. These datasets can be reduced, without lose of information, to Formal Concept Analysis compatible datasets where a binary relation determines if an object has a property or attribute. In this area we can use different methods and properties that have been developed and can not be used in the original dataset.

The mining of negative attributes from datasets has been studied in the last decade to obtain additional and useful information. There exists an exhaustive study around the notion of negative association rules between sets of attributes. However, in Formal Concept Analysis, the needed theory for the management of negative attributes is in an incipient stage. We proposed in a previous work an algorithm, based on the NextClosure algorithm [1], that allows obtaining mixed implications. The proposed algorithm returns a feasible and complete mixed implicational system by performing a reduced number of requests to the formal context.

Knowledge discovering is nowadays a well established discipline focussed on the development of tools and techniques to reveal useful information hidden in data sets. Its main goal is to detect patterns to improve decision making and is approached using pattern recognition, clustering, association and classification. Part of these patterns is expressed as implications (or association rules) which allow us to address information using a formal notation and to manage them syntactically by using logic.

In Formal Concept Analysis we have 3 related items such that we can manage them to obtain knowledge: Implications, Concepts and Context.

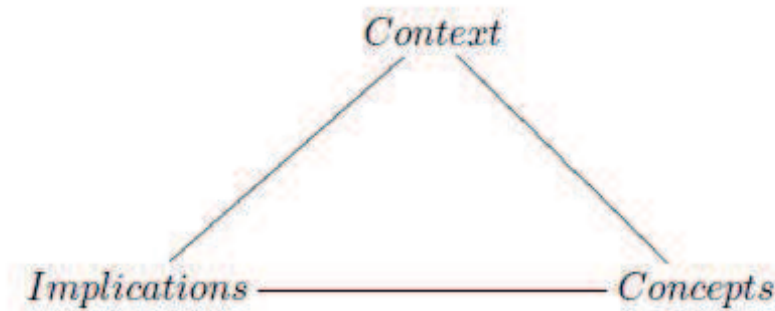

Figure 1. Items in Formal Concept Analysis

Contexts are data sets with a binary relation **I**, between a set of objects **G** and their attributes **M**.

Implications are formulas in the form A→B where A and B are subsets of a certain set (universe) of attributes **M**. Both subsets of attributes are considered to be conjunctive cubes, i.e. $A = a_1 \wedge \ldots \wedge a_n$ and $B = b_1 \wedge \ldots \wedge b_m$.

Concepts are closed sets that have a pair formed by a subset of objects and a subset of attributes that are related with two mappings that are called derivation operators. We can calculate the closure of a subset of attributes with respect to a set of implications and establish the concept associated.

Formal Concept Analysis have a extended set of applications in real life that allows to manage information to extract knowledge from a context, a set of implications or a set of concepts: Social Networks, Marketing, Recommendation Systems,...

We are going to focus in the existence of logics developed to specify and manage sets of implications. The pioneer of these logics was the one introduced by W.W. Armstrong [2], which was proven to be sound and complete.

In this work we use implications in the area of formal concept analysis [3] and assume the common interpretation in this environment: given a formal context **K** over a set of attributes Ω, the implication A→B asserts that any object which have all the A attributes, also has all the B attributes. Although we focus on knowledge discovery techniques to extract implications in formal concept analysis, this problem is similar to the extraction of functional dependencies or association rules from an arbitrary data set.

One of the former researchers who points out the importance of this problem in datasets was H. Mannila [4] and it was also studied by other researches from the database areas like S. Navathe [5]. In this work we will address the mining of implications with negation. The extended implications allow us to relate items which conflict with each other. While classical implications express that {``*cyclist with short and sharp accelerations are great sprinters*''}, implications with negations allow us to express that {``*cyclist with short and sharp accelerations are not great climbers*''}.

Since implication formulas are built using a binary connective which relates two conjunctive clauses, negation is only considered at the attribute level, i.e. the negation $\overline{A}$ is considered the conjunctive clause of its negated attributes: $\overline{a}_1 \wedge \ldots \wedge \overline{a}_n$. We want to increase the expressiveness to exploit all the properties that Formal Concept Analysis had, but not result in much higher costs.

Notice that extended implications cannot be considered as the negation of a classical implication. We do not want to express that a certain implication does not holds but the evidence of the absence of a certain attribute because this is a solved question.

Applications of mixed attributes are usual in data mining but not in Formal Concept Analysis. There are some usual approaches, for example detecting errors [6], which provide benefits to knowledge discovery.

## 2. PRELIMINARIES

In this section, the basic notions related with Formal Concept Analysis (FCA) [7] and attribute implications are briefly presented. See [3] for a more detailed explanation.

A *formal context* is a triple **K**=<G,M,I> where **G** and **M** are finite non-empty sets and I⊆GxM is a binary relation. The elements in G are named objects, the elements in M attributes and <g,m>∈ I means that the object g has the attribute m. From this triple, two mappings ↑:$2^G$→ $2^M$ and ↓:$2^M$→$2^G$, named derivation operators, are defined as follows: for any A⊆G and B⊆M,

$A^\uparrow$ = {m∈ M| for each g∈ A: <g, m>∈ I}

$B^\downarrow$ = {g∈ G| for each m∈ B: <g, m>∈ I}

$A^\uparrow$ is the subset of all attributes shared by all the objects in A and $B^\downarrow$ is the subset of all objects that have the attributes in B. The pair (↑,↓) constitutes a Galois connection between $2^G$ and $2^M$ and, therefore, both compositions are closure operators.

A pair of subsets <A, B> with A⊆G and B⊆M such $A^\uparrow$=B and $B^\downarrow$=A is named a *formal concept*. A is named the *extent* and B the *intent* of the concept. These extents and intents coincide with closed sets wrt the closure operators because $A^{\uparrow\downarrow}$=A and $B^{\downarrow\uparrow}$=B. Thus, the set of all the formal concepts have a lattice structure, named *concept lattice*}, with the relation

<$A_1$, $B_1$> ≤ <$A_2$, $B_2$> if and only if $A_1$⊆$A_2$ (or equivalently, $B_2$⊆$B_1$)

The concept lattice can be characterized in terms of attribute implications. An *attribute implication* is an expression A→B where A, B ⊆ M and it holds in a formal context if $A^{\downarrow} \subseteq B^{\downarrow}$. That is, any object that has all the attributes in A has also all the attributes in B. It is well known that the sets of attribute implications that are satisfied by a context satisfy the Armstrong's Axioms:

[Ref] Reflexivity: If B ⊆ A then ⊢ A→B.

[Augm] Augmentation: A→B ⊢ A∪C→B∪C.

[Trans] Transitivity: A→B, B→C ⊢ A→C.

A set of implications $\mathbb{B}$ is an *implicational system* for K if: (1) any implication from $\mathbb{B}$ holds in K and (2) any implication that K satisfies follows (can be inferred) by using Armstrong's Axioms from $\mathbb{B}$.

## 3. OUR APPROACH FOR EXTENDING FORMAL CONCEPT ANALYSIS WITH NEGATIVES ATTRIBUTES

From now on, the set of all the attributes is denoted by M and its elements by the letter m possibly with sub indexes. The elements in $M \cup \overline{M}$ are going to be denoted by the first letters in the alphabet: a, b, c … So, the symbols a, b, c … could represent positive or negative attributes. Capital letters A, B, C … denote subsets of $M \cup \overline{M}$. If A ⊆ $M \cup \overline{M}$, then $\overline{A}$ denotes the set of the opposite of attributes in A. That is, $\overline{A} = \{\overline{a} \mid a \in A\}$ where $\overline{\overline{a}} = a$. Moreover, for A ⊆ $M \cup \overline{M}$, the following sets are defined:

Pos(A) = {m ∈ M | m ∈ A}

Neg(A) = {m ∈ M | $\overline{m}$ ∈ A}

Tot(A) = Pos(A) ∪ Neg(A)

and, therefore, Pos(A), Neg(A), Tot(A) ⊆ M.

The traditional derivation operators defined in Formal Concept Analysis are modified in [8] to consider the new framework.

**Definition 1:** *Mixed Concept-forming Operators*. Let K=<G, M, I> be a formal context. We define the operators $\Uparrow: 2^G \to 2^{M \cup \overline{M}}$ and $\Downarrow: 2^{M \cup \overline{M}} \to 2^G$ as follows: for A ⊆ G and B ⊆ $M \cup \overline{M}$,

$A^{\Uparrow} = \{m \in M | <g, m> \in I$ for all $g \in A\} \cup \{\overline{m} \in \overline{M} | <g, m> \notin I$ for all $g \in A\}$

$B^{\Downarrow} = \{g \in G | <g, m> \in I$ for all $m \in B\} \cap \{g \in G | <g, m> \notin I$ for all $\overline{m} \in B\}$

We prove in [9] that the adaptation of the derivation operators form a Galois connection.

**Theorem 1**: Let K=<G, M, I> be a formal context. The pair of derivation operators ($\Uparrow, \Downarrow$) introduced in Definition 1 is a Galois Connection.

**Example 1:** Chemical analysis of oil in motors allows detecting problems inside them.

Table 1. Presence of chemical residuum in analysis of oil.

|   | Fe | Cu | Pb | Al | Si |
|---|----|----|----|----|----|
| $o_1$ | 0 | 1 | 1 | 1 | 0 |
| $o_2$ | 0 | 1 | 0 | 1 | 0 |
| $o_3$ | 1 | 0 | 0 | 1 | 0 |
| $o_4$ | 0 | 0 | 0 | 0 | 1 |

$\{o_1,o_2\}^\uparrow =\{Cu, Al\}\neq\{\overline{Fe}, Cu, Al, \overline{Si}\}=\{o_1,o_2\}^\Uparrow$

Also, we adapt the definitions of formal concept and implication to allow the use of negative attributes inside them.

**Definition 2**: Let K=<G, M, I> be a formal context. A *mixed formal concept* in K is a pair of subsets <A,B> with $A \subseteq G$ and $B \subseteq M \cup \overline{M}$ such $A^\Uparrow$=B and $B^\Downarrow$=A.

**Definition 3**: Let K=<G, M, I> be a formal context and let $A,B \subseteq M \cup \overline{M}$, the context K satisfies a *mixed attribute implication* A→B, if $A^\Downarrow \subseteq B^\Downarrow$.

There are special cases of mixed attributes subsets that have important properties for the proposed algorithm and reduce the complexity of it. Missaoui [10] works with intents that have empty support, i.e. have any m∈ M such that m$\overline{m}$ in the intent, so these subsets of attributes doesn't appear in the real life. We want to discard them because they are not useful for our purposes.

**Definition 4**: Let K=<G, M, I> be a formal context and a set $A \subseteq M \cup \overline{M}$ is named consistent set if Pos(A)∩Neg(A)=∅. The set of all consistent sets is denoted **Ctts**.

**Definition 5**: Let K=<G, M, I> be a formal context and a set $A \subseteq M \cup \overline{M}$ is said to be full consistent set if A∈ Ctts and Tot(A)=M.

We are going to extend the logic that propose Armstrong's Axioms with a sound and complete system proposed by Missaoui with a global framework and not for particular problems.

Therefore, the following axioms are added to the Armstrong's axioms: for all a, b∈ $M \cup \overline{M}$ and $A \subseteq M \cup \overline{M}$,

[Cont] Contradiction: $\vdash a\overline{a} \rightarrow M\overline{M}$.

[Rft] Reflection: Aa→ b $\vdash A\overline{b} \rightarrow \overline{a}$.

The algorithm for calculating the implicational system that we propose in [8] uses the axiomatic system and, considering the full consistent sets, the dataset have to be checked $3^{|M|}$-$2^{|M|}$ times in the worst case.

```
Function Closed(A, 𝔅): boolean
    Data: A ∈ ℭtts, and 𝔅 being a set of mixed implications.
    Result: 'true' if A is closed wrt 𝔅 or 'false' otherwise.
1   begin
2       foreach B → C ∈ 𝔅 do
3           if B ⊆ A and C ⊄ A then
4               exit and return false
5           if B \ A = {a}, A ∩ C̄ ≠ ∅, and ā ∉ A then
6               exit and return false
7       return true
8   end
```

```
Algorithm 2: Mixed Implications Mining
    Data: 𝕂 = (G, M, I)
    Result: Σ set of implications
1   begin
2       Σ := ∅;
3       Y := ∅;
4       while Y < M do
5           foreach X ⊆ Y do
6               A := (Y \ X) ∪ X̄;
7               if Closed(A, Σ) then
8                   C := A^{⇓⇑};
9                   if A ≠ C then
10                      Σ := Σ ∪ {A → C \ A}
11          Y := Next(Y) // i.e. successor of Y in the lectic order
12      return Σ
13  end
```

## 4. ALGORITHM FOR MIXED CONCEPTS AND MIXED ATTRIBUTE IMPLICATIONS

If we want to calculate both, mixed implicational system and set of mixed concepts, we can extend algorithm 2 [8]. We know that we only have to modify a few lines because it works with closures.

To adapt the algorithm, we have to pay attention to the condition *While Y< M* because it has to be changed to *While Y≠M* for detecting the concepts, so the time of calculating arise. Because we only have to add full consistent sets, that are the minimal concepts that have the intents of the objects, we don't need to change this condition, and only add these intents.

If we choose the first option we have to do $2^{|M|}$ operations that always be bigger or equal than |G|. So, the modifications are in lines 12-13, that take the closed sets that we ignore in Algorithm 1, and lines 15-16 that add to the sets of intents of concepts the full consistent sets.

```
Algorithm 3: Mixed Implications and Concepts Mining
    Data: K = ⟨G, M, I⟩
    Result: Σ set of implications, B set of intents
 1  begin
 2      Σ := ∅;
 3      B := ∅;
 4      Y := ∅;
 5      while Y < M do
 6          foreach X ⊆ Y do
 7              A := (Y ∖ X) ∪ X̄;
 8              if Closed(A, Σ) then
 9                  C := A^⇓⇑;
10                  if A ≠ C then
11                      Σ := Σ ∪ {A → C ∖ A}
12                  else
13                      B := B ∪ C
14          Y := Next(Y) // i.e. successor of Y in the lectic order
15      foreach g ∈ G do
16          B := B ∪ g^⇑
17      return Σ, B
18  end
```

With this algorithm we extract all possible implications and concepts from the context. In the worst case, we have to check $3^{|M|}-2^{|M|}+|G|$ times the dataset that not depend of different amount of objects. The dataset could be mixed-clarified and the number of possible objects has a maximum in $2^{|M|}$.

## 5. OPTIMIZING ENGINE DESIGN

In engine simulators, each component has specific measures that determine how these components work. Width of a simple pipe could increase or decrease the amount of fuel that the engine receive, so all this measures are important in the development of each engine.

Each measure could be considered as an attribute of the engine and also, the relation among them determine the global efficiency. All these measures are collected in a dataset that show how the engine works but it is difficult to decide specifically which of them we can change and how we have to change it because we don't know how these changes affect engine performance.

Taguchi [11] refers that development engineers have to apply parameter design methods to make the basic functions approach the ideal functions under real conditions. These design activities should be conducted by research and development departments before the product is finally created. Taguchi's work is based in statistical methods that actually are replaced with neural networks processes or Multi-objective Covariance Matrix Adaptation Evolutionary Strategy methods [12, 13]

The method that we propose tries to focus on these specific parts of the engine that could increase the performance and change these values for obtaining the best result. Researchers could observe that is similar to a self-organized map [14], with the difference that use properties of mixed concept lattices and we only have to change a reduced number of attributes, reducing the cost of the procedure.

Our starting point is different experimental configurations of the engine that experts consider as initial results. Let G be the set of different collected configurations (objects) and M the set of elements that could be measured adding a value called *objective function* (attributes). Value of the objective function depends on different values of attributes, for example number of revolutions, but could not be associated with a particular equation. However, in the next example, for demonstration purposes, we have to consider a polynomial function.

**Example 2**: Let consider domain of possible values [0, 10] for all attributes and the objective function is $a^2+(b-5)+(c-5)-(d-5)-e^2$, so objective domain is [-115,115]. The attributes with higher weights are a and e, so changing these values we can modify the objective function in different ways.

Let's consider that experts determine that the optimal value for the objective is 100 and we want to know what are the preferred size for different parts of the engine for obtaining the nearest value to the proposed optimal.

Table 2: Dataset with engine configurations

| T | Part a | Part b | Part c | Part d | Part e | Objective |
|---|---|---|---|---|---|---|
| $c_1$ | 3 | 1 | 7 | 2 | 3 | 1 |
| $c_2$ | 4 | 3 | 2 | 3 | 4 | -3 |
| $c_3$ | 6 | 5 | 4 | 4 | 2 | 32 |
| $c_4$ | 6 | 0 | 3 | 6 | 5 | 3 |
| $c_5$ | 2 | 7 | 7 | 3 | 1 | 9 |

Let T be a dataset with size $|G|*|M|$ that considers configurations in rows and measures of each element of the engine in columns. Elements in T, ($t_{x,y}$ with $x \in G$, $y \in M$) are numbers that represents each different value for the attributes but not necessary represents a formal context.

We have to point to the objective function and determine the selected value for this function that experts consider as optimal. Different values for objective function could be near from the optimal value but we need the optimal or, in cases where it is not possible, the nearest one.

A specific number k, with $1<k<|G|$, have to be selected that determine the size of the control fixed group. This group $G_k$ is compounded by k configurations which have the nearest objective values with respect to the optimal value. For each attribute $m \in M$, we can fix an interval $In_m=[l_m,u_m]$ being $l_m=\inf(t_{x,m})$ and $u_m=\sup(t_{x,m})$ with $x \in G_k$.

Now we can adapt T to a compatible dataset $T^1$ in Formal Concept Analysis in this way. Let define the values of the adapted dataset as $t^1_{x,y}=1$ if $t_{x,y} \in In_y=[l_y,u_y]$; 0 otherwise

This adaptation represents all the values that are in the fixed intervals with a binary relation. It is obvious that, for all $g \in G_k$, $m \in M$, $t^1_{g,m}=1$ so we have to focus in the absence of attributes in the sense of Formal Concept Analysis. The absence of attributes means that these values are out of the fixed intervals (fixed by $G_k$) and points to possible candidates for replacing their values.

**Example 3**: Fixing k=2 and optimal value for objective 100 we have that $G_k=\{c_3,c_5\}$ and $In_a=[2,6]$, $In_b=[5,7]$, $In_c=[4,7]$, $In_d=[3,4]$, $In_e=[1,2]$

Table 3: Adapted formal context to T

| $T^1$ | Part a | Part b | Part c | Part d | Part e | Objective |
|---|---|---|---|---|---|---|
| $c_1$ | 1 | 0 | 1 | 0 | 0 | 1 |
| $c_2$ | 1 | 0 | 0 | 1 | 0 | -3 |
| $c_3$ | 1 | 1 | 1 | 1 | 1 | 32 |
| $c_4$ | 1 | 0 | 0 | 0 | 0 | 3 |
| $c_5$ | 1 | 1 | 1 | 1 | 1 | 9 |

The problem that we have to solve is what $t_{x,y} \in T$ we have to change and what values could have $t_{x,y}$.

Using Algorithm 3 with the context $<G-G_k, M, T^1>$, we can obtain the mixed concepts lattice. The top element is the mixed concept $<G-G_k, (G-G_k)^{\Uparrow}>$. Since the lattice has an order, the upper concept that has the less number of negative attributes (not zero) determines the selected set of attributes to change. We are going to call this subset of attributes N.

For avoiding local fix points for the objective function, we have to consider 4 different values for each attribute that we are going to change in T: original value, average between minimum possible value for attribute y and $l_y$, average between maximum possible value for attribute y and $u_y$, and average for the interval $In_y$, for each $y \in N$, i.e. for each configuration, we have to consider $4^{|N|}$ combinations.

If changing attributes in N we do not obtain better objective values, we have to follow in order the mixed concept lattice and get the next mixed concept for obtaining the new N. If proposed changes do not offer a better objective value for all the objects that are not in $G_k$, we have to finish our search.

This method have to be repeated recursively while distance between objective value and optimal value in $G_k$ are higher than distance between objective value and optimal value in $G-G_k$ after checking all possible combinations.

**Example 4**: Obtaining mixed concept lattice for $G-G_k$, the top element of the lattice point to the concept with attributes $a\overline{be}$. Parameter a is positive so it is not necessary to be changed. Focusing in negative attributes, we have to change parameters b, e in $c_1$, $c_2$ and $c_4$, resulting Table 4.

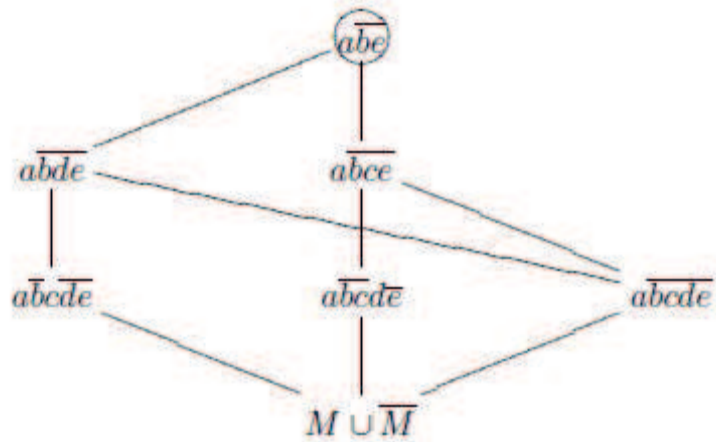

Figure 2. Mixed concept lattice for $T^1$ restricted to $c_1$, $c_2$ and $c_4$

Table 4. Adapted formal context for T

| T | Part a | Part b | Part c | Part d | Part e | Pre. Obj. | New Obj. |
|---|---|---|---|---|---|---|---|
| $c_1$ | 3 | 8.5 | 7 | 2 | 0.5 | 1 | 17.25 |
| $c_2$ | 4 | 8.5 | 2 | 3 | 0.5 | -3 | 18.25 |
| $c_3$ | 6 | 5 | 4 | 4 | 2 | 32 | |
| $c_4$ | 6 | 8.5 | 3 | 6 | 0.5 | 3 | 36.25 |
| $c_5$ | 2 | 7 | 7 | 3 | 1 | 9 | |

We can observe that all the objectives values have changed and, in this case, $c_1$, $c_2$ and $c_4$ have an objective value higher than $c_5$ that was in $G_k$ in the first step, and value of $c_4$ is higher than $c_3$ so the distance to the optimal value decreases. Values for b increases and values for e decreases.

In the next step, $G_k=c_3,c_4$ and, recursively, we change in second step groups a, c, d; third step a, c, d; fourth step a, c, d; fifth step a, b, d, e and we find that has the same objective value that optimal selected by experts.

In 2 steps we have to check $3*(4^2)+9*(4^3)+4^4=880$ combinations that are less than all possible combinations for each configuration, for example, if we only take into account positive integers, there are $10^5=100000$ combinations. If we consider decimals our method is not affected but the possible combinations grow.

In this example, if we consider minimum value for y instead of average between minimum possible value for y and $l_y$ and maximum value for y instead of average between maximum possible value for y and $u_y$, for each $y \in N$, we only have to do 112 combinations to reach the proposed optimal value, but this example have a particular solution using interval extremes, that is not the general case.

This method could be included in simulators as an artificial intelligence tool for adjusting parameters. Being an automatic tool, engineers do not need to examine and calibrate the simulator trying to detect what parameter has to change and what value is needed. Configuration for optimal value is not always reached due to parameters restrictions and the proposed original configurations, but a close value could be obtained.

As other methods, detecting local fix points is a problem that can be solved with a different set of initial configurations as a start point.

## 6. EXPERIMENTS

Due to legal problems in the use of confidential data, we can not check the method proposed in the previous section in real simulators, so we have to do simulations with different equations. These equations have been chosen for its simplicity of calculation and complexity and diversity of results according to the parameters entered.

### 6.1. Experiment 1

Table 5. Results from Experiment 1

| k | |G| | Objective | Init. Dist. | Final Dist. | Iter. | % Reduc. |
|---|---|---|---|---|---|---|
| 2 | 5 | 930 | 577.73 | 0.29 | 4.1 | 0.05 |
| 2 | 10 | 930 | 418.70 | 0.75 | 4.7 | 0.18 |
| 2 | 20 | 930 | 197.20 | 1.11 | 4.2 | 0.56 |
| 3 | 5 | 930 | 796.80 | 0.68 | 4.1 | 0.09 |
| 3 | 10 | 930 | 224.65 | 0.30 | 5.6 | 0.14 |
| 3 | 20 | 930 | 117.62 | 0.41 | 3.6 | 0.35 |
| 2 | 5 | 27030 | 14403.09 | 0.05 | 11.0 | 0.00 |
| 2 | 10 | 27030 | 14339.90 | 0.59 | 18.0 | 0.00 |
| 2 | 20 | 27030 | 10332.63 | 0.78 | 19.3 | 0.01 |
| 3 | 5 | 27030 | 15115.54 | 0.10 | 11.0 | 0.00 |
| 3 | 10 | 27030 | 13309.24 | 0.44 | 14.0 | 0.00 |
| 3 | 20 | 27030 | 13819.18 | 0.34 | 21.5 | 0.00 |
| 2 | 5 | 27900 | 19507.26 | 0.12 | 19.8 | 0.00 |
| 2 | 10 | 27900 | 15730.50 | 0.43 | 31.1 | 0.00 |
| 2 | 20 | 27900 | 14026.64 | 0.96 | 46.7 | 0.01 |
| 3 | 5 | 27900 | 16117.79 | 0.05 | 15.5 | 0.00 |
| 3 | 10 | 27900 | 17525.48 | 0.81 | 34.7 | 0.01 |
| 3 | 20 | 27900 | 11296.91 | 0.76 | 38.1 | 0.01 |

The size of M was fixed in 9 attributes which have their domains in [0,10] that allows decimal values. The objective function for this experiment is $\sum_{i=1}^{3} t_{x,i} + (\sum_{i=4}^{6} t_{x,i})^2 + (\sum_{i=7}^{9} t_{x,i})^3$ for each configuration $x \in G$. This function has values in [0,27930] and has 3 groups where we can interexchange values with the same result, i.e. if we change values in $t_{x,4}$ and $t_{x,5}$ objective function does not change.

We are going to use different values for k and |G| for the objectives 930, 27030 and 27900 that are the different configurations obtained decreasing one of the three groups to their lower values and the other groups to the upper values. 10 experiments were done with different random initial values for each fixed value for k, |G| and objective and the average results are presented in Table 5. Initial distance is the difference in absolute distance from the best original objective value to the optimal value, final distance is the absolute distance from the best objective value to the optimal value after applying the method, iterations are the number of recursive executions of the method and reduction is the proportion between initial distance and final distance.

In this experiment we can see that, for the objective 930, there are multiple combinations of values that reach this objective and our method has difficult to locate one of them specifically due to the several fix points. For objectives 27030 and 27900, the number of possible combinations is lower, so our method can focus in particular combinations with a relevant percentage of reduction from the original distance to the optimal value. Is an interesting observation that low size of |G| has the nearest solutions to the optimal value.

### 6.2. Experiment 2

The size of M was fixed in 9 attributes which have their domains in [0,10] that allows decimal values. The objective function for this experiment is $\sum_{i=1}^{9} (t_{x,i} - i)^2$ for each configuration $x \in G$. This function has values in [0,485] and has a unique combination of values that reach the objective value 0. The difficulty to reach the optimal value in 0 is that each attribute have to get a value that is not in the extremes of its domain.

We are going to use different values of k and |G| simulating 10 experiments with different random initial values for each fixed value for k and |G|. The average results are presented in Table 6.

In this experiment, final distance is, more or less, similar in the average value of each experiment. We have a specific case where k=3 and |G|=5 where the solution is the nearest one to the optimal value. Since initial distance is lower when |G| increase, the percentage of reduction grows, but we can see that there is not a special relation between k and |G| in this experiment for reaching the nearest solution to the optimal value.

Table 6. Results from Experiment 2

| k | |G| | Init. Dist. | Final Dist. | Iter. | % Reduc. |
|---|---|---|---|---|---|
| 2 | 5 | 72.87 | 2.68 | 5.8 | 3.66 |
| 2 | 10 | 52.87 | 2.67 | 7.5 | 5.06 |
| 2 | 20 | 54.30 | 2.66 | 8.5 | 5.33 |
| 2 | 30 | 44.40 | 2.70 | 9.8 | 6.09 |
| 3 | 5 | 93.97 | 0.73 | 8.0 | 0.77 |
| 3 | 10 | 66.30 | 2.56 | 8.9 | 3.86 |
| 3 | 20 | 52.78 | 2.36 | 8.1 | 4.48 |
| 3 | 30 | 45.14 | 2.11 | 9.2 | 4.69 |
| 4 | 5 | 89.72 | 2.03 | 8.7 | 2.27 |
| 4 | 10 | 64.41 | 2.56 | 7.8 | 3.97 |
| 4 | 20 | 47.16 | 2.29 | 9.3 | 4.86 |
| 4 | 30 | 47.08 | 2.46 | 10.9 | 5.21 |

## 7. CONCLUSIONS AND FUTURE WORK

In this paper, a new algorithm for obtaining mixed concept lattices from a context is introduced. This algorithm is an important tool that let us to consider how attributes in the formal context are related among them, focusing our attention to relevant information.

The algorithm also produces an implicational system that we do not use in this work but is a complementary tool that gives us specific information about relations between attributes and was the essential idea to design the method that we propose.

This method for adjusting parameters in simulators is not a new idea, but we want to extract the knowledge hidden in datasets using formal concept analysis and, in particular, using the knowledge that negative attributes provide.

Experiments let us to check that this method could be an interesting tool for simulators with a relevant reduction in production and development costs.

In a future work, we want to develop new algorithms to reduce the time of calculating the mixed concepts lattice adapting different algorithms from Formal Concept Analysis.

Values of k and |G| have to be studied with different situations, for example changing |M|, checking what the best initial values for these parameters are.

## ACKNOWLEDGEMENTS

Supported by grant TIN2011-28084 and TIN2014-59471 of the Science and Innovation Ministry of Spain, co-funded by the European Regional Development Fund (ERDF).

We want to thanks Dr. Juan Cabrera and Dr. Juan Castillo for the approach to solving the problem discussed in this article.